\documentclass[10pt]{amsart}
\usepackage{amsmath}
\usepackage{amssymb}
\newcommand{\E}{\mathbb E}
\newcommand{\eR}{\mathbb R}

\newcommand{\N}{\mathbb N}
\newcommand{\rf}[1]{{\rm(\ref{#1})}}
\newtheorem{example}{Example}
\newtheorem{definition}{Definition}
\newtheorem{remark}{Remark}
\newtheorem{theorem}{Theorem}
\newtheorem{lemma}{Lemma}
\newtheorem{proposition}{Proposition}

\begin{document}

\title{Inequalities of the Hermite-Hadamard type involving numerical differentiation formulas}

\author{Andrzej Olbry\'s       
         %etc.
}

\address{Andrzej Olbry\'s  \\University of Silesia
              Institute of Mathematics \\
              Bankowa 14\\
              40-007 Katowice\\
              Poland }             
\email{andrzej.olbrys@wp.pl}           %  \\

\author{Tomasz Szostok} 
             \address{Tomasz Szostok\\ University of Silesia\\ Institute of Mathematics \\
             Bankowa 14\\
              40-007 Katowice\\
              Poland}
              \email{tszostok@math.us.edu.pl}

\begin{abstract}
 We observe that the Hermite-Hadamard inequality written in the form 
$$f\left(\frac{x+y}{2}\right)\leq\frac{F(y)-F(x)}{y-x}\leq\frac{f(x)+f(y)}{2}$$
may be viewed as an inequality between two quadrature operators $f\left(\frac{x+y}{2}\right)$
$\frac{f(x)+f(y)}{2}$ and a differentiation formula $\frac{F(y)-F(x)}{y-x}.$
We extend this inequality, replacing the middle term by more complicated ones. As it turns out
in some cases it suffices to use Ohlin lemma as it was done in a recent paper \cite{Rajba}
however to get more interesting result some more general tool must be used.  
To this end we use Levin-Ste\v{c}kin theorem which provides necessary and sufficient 
conditions under which inequalities of the type we consider are satisfied.

\end{abstract}
\keywords{Hermite-Hadamard inequality, differentiation formulas, convex functions }
\subjclass[2010]{26A51, 26D10, 39B62}
\maketitle

\section{Introduction}
\label{intro}
We shall obtain some class of  inequalities of the Hermite-Hadamard type. 
First we write the classical Hermite-Hadamard inequality    
\begin{equation}
\label{HH}
f\left(\frac{x+y}{2}\right)\leq\frac{1}{y-x}\int_{x}^yf(t)dt\leq\frac{f(x)+f(y)}{2}
\end{equation}
(see \cite{DP} for many generalizations and applications of \rf{HH}).
Now, let us write \rf{HH} in the form    
\begin{equation}
\label{HH1}
f\left(\frac{x+y}{2}\right)\leq\frac{F(y)-F(x)}{y-x}\leq\frac{f(x)+f(y)}{2}.
\end{equation}
This inequality is, clearly, satisfied by every convex function $f$ and its primitive function $F.$
However  \rf{HH1} may be viewed as an inequality involving two types of expressions used, respectively, in numerical integration and differentiation. Namely 
$f\left(\frac{x+y}{2}\right)$ and $\frac{f(x)+f(y)}{2}$ are the simplest quadrature formulas 
used to approximate the definite integral, whereas 
$\frac{F(y)-F(x)}{y-x}$ is the simplest expression used to approximate the derivative of $F.$
Moreover, as it is known from numerical analysis, if $F'=f$ then the following equality is satisfied
\begin{equation}
\label{er1}
f(x)=\frac{F(x+h)-F(x-h)}{2h}-\frac{h^2}{6}f''(\xi)
\end{equation}
for some $\xi\in(x-h,x+h).$ This means that \rf{er1} provides an alternate proof of \rf{HH}
(for twice differentiable $f$).

This new formulation of the Hermite-Hadamard inequality inspires us to replace
the middle term of Hermite-Hadamard inequality by more complicated expressions than those used in \rf{HH}.    
In the paper \cite{Szostok} all numbers $a,\alpha,\beta\in[0,1]$ such that for all convex functions 
$f$ the inequality 
$$af(\alpha x+(1-\alpha )y)+(1-a)f(\beta x+(1-\beta y)\leq \frac{1}{y-x}\int_{x}^yf(t)dt$$
is satisfied and all $a,b,c,\alpha\in(0,1)$ with $a+b+c=1$ for which we have
$$bf(x)+cf(\alpha x+(1-\alpha)y)+df(y)\geq\frac{1}{y-x}\int_{x}^yf(t)dt.$$ This means that 
the expressions   $f\left(\frac{x+y}{2}\right)$ and $\frac{f(x)+f(y)}{2}$ were replaced by longer ones while the integral mean remained unchanged. In the current paper
we do the opposite job, namely, we are going to prove inequalities of the type 
$$f\left(\frac{x+y}{2}\right)\leq\frac{a_1F(x)+a_2F(\alpha x+(1-\alpha)y)+a_3F(\beta x+(1-\beta)y)+a_4F(y)}{y-x}$$
and
$$\frac{a_1F(x)+a_2F(\alpha x+(1-\alpha)y)+a_3F(\beta x+(1-\beta)y)+a_4F(y)}{y-x}\leq\frac{f(x)+f(y)}{2}$$
where $f:[x,y]\to\eR$ is a convex function, $F'=f,$ $\alpha,\beta\in(0,1)$ and $a_1+a_2+a_3+a_4=0.$  

\section{Preliminaries}

In recent papers \cite{Rajba} and \cite{Szostok}  Ohlin lemma on convex stochastic ordering was used to obtain inequalities of the Hermite-Hadamard type. In the current paper
we are going to work in a similar spirit. Thus now we cite this lemma.
\begin{lemma} (Ohlin \cite{Ohlin})
Let $X_1,X_2$ be two random variables such that $\E X_1=\E X_2$  and let  $F_1,F_2$ be their cumulative distribution functions.
If $F_1,F_2$ satisfy for some $x_0$ the following inequalities
\begin{equation}
\label{x0}
F_1(x)\leq F_2(x)\;{\rm if}\;x<x_0\;\;{\rm and}\;\;F_1(x)\geq F_2(x)\;{\rm if}\;x>x_0 
\end{equation}
then 
\begin{equation}
\label{m}
\E f(X_1)\leq \E f(X_2)
\end{equation}
 for all continuous and convex functions $f:\eR\to\eR.$
 \label{OL}
\end{lemma}
However, in the present approach we are going to use a result from \cite{LS}, (see also \cite{NP}
Theorem 4.2.7).

\begin{theorem}(Levin, Ste\v{c}kin)
\label{LS}
Let $F_1,F_2:[a,b]\to\eR$ be two functions with bounded variation such that $F_1(a)=F_2(a).$
Then, in order that 
$$\int_a^b f(x)dF_1(x)\leq \int_a^b f(x)dF_2(x),$$
for all continuous and convex $f$, 
it is necessary and sufficient that $F_1$ and $F_2$ verify the following three conditions:
\begin{equation}F_1(b)=F_2(b),\end{equation}
\begin{equation}\int_a^x F_1(t)dt\leq\int_a^x F_2(t)dt,\;x\in(a,b),
\label{FGi}
\end{equation}
\begin{equation}\int_a^b F_1(t)dt=\int_a^b F_2(t)dt.
\label{FGe}
\end{equation}
\end{theorem}

\begin{remark}\rm Observe that if measures $\mu_1,\mu_2$ corresponding to the random variables
occurring in Ohlin's lemma are concentrated on the interval $[x,y]$ then Ohlin's lemma  is an easy consequence of Theorem \ref{LS}. Indeed, $\mu_1,\mu_2$ are probabilistic measures thus we have
$F_1(x)=F_2(x)=0$ and $F_1(y)=F_2(y)=1.$ Moreover $\E X_1=\E X_2$ yields \rf{FGe} and, from the inequalities \rf{x0}, we get \rf{FGi}.  
 \end{remark}

Now we shall use Theorem \ref{LS} to make an observation which is more general than 
Ohlin lemma and concerns the situation when functions $F_1,F_2$ have more crossing points than one.
First we need the following definition.
\begin{definition}
Let $F_1,F_2:[a,b]\to\eR$ be functions and let $a=x_0<x_1<\dots<x_n<x_{n+1}=b.$ We say that the pair $(F_1,F_2)$ crosses $n-$times (at points $x_1,\dots,x_n$) if the inequalities 
\begin{equation}
\label{xi}
F_1(x)\leq F_2(x),\;x\in(x_i,x_{i+1})\;{\rm and}\;F_1(x)\geq F_2(x),\;x\in(x_{i+1},x_{i+2})
\end{equation}
where $i$ is even, are satisfied and
$$\int_{x_i}^{x_{i+1}}[F_1(x)-F_2(x)]dx\neq 0,\qquad i=0,1,...,n.$$
\end{definition}

\begin{lemma}
\label{mainl}
Let $F_1,F_2:[a,b]\to\eR$ be two functions with bounded variation such that $F_1(a)=F_2(a)$ and $F_1(b)=F_2(b)$
let $x_1,\dots,x_n\in(a,b)$ and let $(F_1,F_2)$ cross at $x_1,\dots,x_n\in(a,b).$

(i) If $n$ is even then the inequality
\begin{equation}
\label{maini}
\int_a^b f(x)dF_1(x)\leq \int_a^b f(x)dF_2(x)
\end{equation}
  is not satisfied by all continuous and convex $f.$
  
(ii) If $n$ is odd then we define numbers $A_i$ by the following formulas
$$A_i:=\int_{x_i}^{x_{i+1}}[F_{1}(x)-F_{2}(x)]dx$$ for odd $i$ and 
$$A_i:=\int_{x_i}^{x_{i+1}}[F_{2}(x)-F_{1}(x)]dx$$ for even $i.$
Inequality \rf{maini} is satisfied for all continuous and convex $f$  if and only if the following
inequalities hold true
$$A_0\geq A_1,$$
$$A_0-A_1+A_2\geq A_3$$
$$A_0-A_1+A_2-A_3+A_4\geq A_5$$
$$\vdots$$
$$A_0-A_1+A_2-A_3+A_4-A_5+\dots+A_{n-4}+A_{n-3}\geq A_{n-2}.$$
\end{lemma} 
{\sc Proof.} First we prove (i). For an indirect proof assume that $n$ is even and that 
\rf{maini} is satisfied for all convex $f.$ Then from  Theorem \ref{LS} (and from the 
definition of $A_0$) we get 
$$A_0=\int_{a}^{x_1}[F_{1}(t)-F_{2}(t)]dt>0.$$ Further, since there is an even number of sign changes, we have
$A_n>0.$ From \rf{FGe} we get 
$$A_0-A_{1}+\cdots+A_n=\int_{a}^b[F_1(t)-F_2(t)]dt=0$$
thus we have
$$\int_{a}^{x_n}[F_{1}(t)-F_{2}(t)]dt=A_0-A_{1}+\cdots-A_{n-1}=-A_n<0$$
which is a contradiction with \rf{FGi}.

Now we shall show the condition (ii) to this end assume that $n$ is odd and that 
\rf{maini} is satisfied for all convex $f.$ Then 
$$A_0-A_1=\int_{a}^{x_2}[F_{1}(t)-F_{2}(t)]dt\geq 0$$
$$A_0-A_1+A_2-A_3=\int_{a}^{x_4}[F_{1}(t)-F_{2}(t)]dt\geq 0$$
$$\vdots$$
$$A_0-A_1+A_2-A_3+A_4+\dots+A_{n-3}-A_{n-2}=\int_{a}^{x_{n-1}}[F_{1}(t)-F_{2}(t)]dt\geq0.$$
Assume, on the other hand that all the above inequalities hold true and take 
$z\in[a,b).$ Then  $z\in [x_{i_{0}},x_{i_{0}+1})$\ for some $i_{0}\in \{0,...,n-1\}$.
If $i_0$ is odd then we have
\begin{eqnarray*}
\int_{a}^{z}[F_{1}(t)-F_{2}(t)]dt=&\sum_{j=0}^{i_{0}}(-1)^j A_{j}+\int_{x_{i_{0}}}^{z}[F_{1}(t)-F_{2}(t)]dt\\
\\
&\geq \sum_{j=0}^{i_{0}+1}(-1)^jA_{j}\geq0.
\end{eqnarray*}
Further if $i_0$ is even then
\begin{eqnarray*}
\int_{a}^{z}[F_{1}(t)-F_{2}(t)]dt=&\sum_{j=0}^{i_{0}}(-1)^j A_{j}+\int_{x_{i_{0}}}^{z}[F_{1}(t)-F_{2}(t)]dt\\
\\
&\geq \sum_{j=0}^{i_{0}}(-1)^jA_{j}\geq0.
\end{eqnarray*}
To finish the proof it is enough to apply the Levin-Ste\v{c}kin theorem.
\section{Results and applications}
In this part of the paper we shall work with expressions of the type 
$$\frac{\sum_{i=1}^{n}a_iF(\alpha_i x+(1-\alpha_i)y)}{y-x}$$
where $\sum_{i=1}^na_i=0.$ Therefore, now we make the following observation,
$l_1$ stands here for the one-dimensional Lebesgue measure.
\begin{proposition}
Let $n\in\N,$ let $\alpha_i\in(0,1)$ satisfy $\alpha_1>\alpha_2>\cdots>\alpha_n,$ let $a_1+a_2+\cdots+a_n=0$ and let $F$ be a differentiable function with $F'=f.$ Then 
$$\frac{\sum_{i=1}^{n}a_iF(\alpha_i x+(1-\alpha_i)y)}{y-x}=\int fd\mu$$
where 
$$\mu(A)=-\frac{1}{y-x}\sum_{i=1}^{n-1}(a_1+\cdots+a_i)l_1(A\cap[\alpha_{i} x+(1-\alpha_{i})y,\alpha_{i+1} x+(1-\alpha_{i+1})y]).$$
\label{l1}
\end{proposition}
{\sc Proof.} Let $x_{j}:=\alpha_{j}x+(1-\alpha_{j})y,\ j=1,...,n$,\ and denote by
\begin{displaymath}
\delta_{k,n}=\left \{ \begin{array}{ll}
1\ \  k\leq n, \\
0\ \  k>n.
 \end{array} \right.
\end{displaymath}
Then, by the definition of $\mu,$ we obtain
$$
\setlength{\arraycolsep}{1pt}
\begin{array}{ll}
\int fd\mu\ & = \sum_{j=1}^{n-1}\int_{x_{j}}^{x_{j+1}}f(t)\mu(dt)=
\sum_{j-1}^{n-1}\Big(-\sum_{k=1}^{j}a_{k} \Big) \frac{F(x_{j+1})-F(x_{j})}{y-x}\\
\\
&=-\sum_{j=1}^{n-1}\Big(\sum_{k=1}^{n-1}\delta_{kj} a_{k} \frac{F(x_{j+1})-F(x_{j})}{y-x}\Big)\\
\\
&=-\sum_{k=1}^{n-1}a_{k}\sum_{j=1}^{n-1}\delta_{kj}  \frac{F(x_{j+1})-F(x_{j})}{y-x}\\
\\
&=-\frac{1}{y-x}\sum_{k=1}^{n-1}a_{k}\sum_{j=k}^{n-1}[F(x_{j+1})-F(x_{j})]\\
\\
&=-\frac{1}{y-x}\sum_{k=1}^{n-1}a_{k}[F(x_{n})-F(x_{k})]\\
\\
&=-\frac{1}{y-x}\Big[\sum_{k=1}^{n-1}a_{k}F(x_{n})-\sum_{k=1}^{n-1}a_{k}F(x_{k})\Big]=\frac{1}{y-x}\sum_{k=1}^{n}a_{k}F(x_{k})
\end{array}
$$

\begin{remark}Taking $F_1(t):=\mu((-\infty,t])$ with $\mu$ from Proposition $1$ we can see that 
\label{R1} 
\begin{equation}
\label{Fint}
\frac{\sum_{i=1}^{n}a_iF(\alpha_i x+(1-\alpha_i)y)}{y-x}=\int fdF_1.
\end{equation}
\end{remark}

Next proposition will show that, in order to get some inequalities of the
Hermite-Hadamard type, we have to use sums containing more than three summands.
\begin{proposition}There are no numbers $\alpha_i,a_i\in\eR,i=1,2,3$ satisfying $1=\alpha_1>\alpha_2>\alpha_3=0$ such that any of the inequalities
$$f\left(\frac{x+y}{2}\right)\leq\frac{\sum_{i=1}^{3}a_iF(\alpha_i x+(1-\alpha_i)y)}{y-x}$$
or 
$$\frac{\sum_{i=1}^{3}a_iF(\alpha_i x+(1-\alpha_i)y)}{y-x}\leq\frac{f(x)+f(y)}{2}$$
is fulfilled by every continuous and convex function $f$ and its antiderivative $F.$
\end{proposition}
{\sc Proof.} Using Proposition \ref{l1}, we can see that 
$$\frac{\sum_{i=1}^{3}a_iF(\alpha_i x+(1-\alpha_i)y)}{y-x}=\int_x^yfd\mu$$
with $$\mu(A)=-\frac{1}{y-x}\bigl(a_1l_1(A\cap[x,\alpha_2x+(1-\alpha_2)y])+$$$$
(a_2+a_1)l_1(A\cap[\alpha_2x+(1-\alpha_2)y,y])\bigr),$$ 
and
$$\frac{\sum_{i=1}^{3}a_iF(\alpha_i x+(1-\alpha_i)y)}{y-x}=\int_x^yf(t)dF_1(t)$$
where
\begin{equation}
\label{F1}
F_1(t)=\mu\{(-\infty,t]\}.
\end{equation}
 Now, let 
$$F_2(t)=\frac{1}{y-x}l_1\{(-\infty,t]\cap[x,y]\}$$ then
$F_1$ lies strictly above or below $F_2$ (on $[x,y]$). This means that 
\begin{equation}
\label{2}
\int_{x}^yF_2(t)dt\neq\int_{x}^yF_1(t)dt.
\end{equation}
But, on the other hand, if 
\begin{equation}
\label{F3}
F_3(t):=\left\{\begin{array}{ll}
0,\ &t<x\\
\frac{1}{2},\ &t\in[x,y)\\
1\ &t\geq y
\end{array}\right.
\end{equation}
and 
\begin{equation}
\label{F4}
F_4(t):=\left\{\begin{array}{ll}
0, & t<\frac{x+y}{2}\\
1, & t\geq\frac{x+y}{2}
\end{array}\right.
\end{equation}
then 
$$\int_{x}^yF_2(t)dt=\int_{x}^yF_3(t)dt=\int_{x}^yF_4(t)dt=\frac{y-x}{2}.$$
This, together with \rf{2}, shows that neither 
$$\int_{x}^yfdF_2\leq\int_{x}^yfdF_3$$ 
nor
$$\int_{x}^yfdF_2\geq\int_{x}^yfdF_4$$ 
is satisfied. To complete the proof it suffices to observe that 
$$\int_{x}^yfdF_3=\frac{f(x)+f(y)}{2}\;{\rm and}\;\int_{x}^yfdF_4=f\left(\frac{x+y}{2}\right).$$
\begin{remark} Observe that the assumptions $\alpha_1=1$ and $\alpha_3=0$ are essential.
For example, it follows from Ohlin lemma that inequality 
$$f\left(\frac{x+y}{2}\right)\leq\frac{-3F(\frac34 x+\frac 14 y)+\frac{25}{11}F(\frac{11}{20}x+\frac{9}{20}y)+\frac{8}{11}F(y)}{y-x}\leq\frac{1}{y-x}\int f(t)dt$$ 
is satisfied by all continuous and convex functions $f$ (where $F'=f$). Clearly there are many 
more examples of inequalities of this type.
\end{remark}

Now we shall present an example showing that it is possible to prove some inequalities of the Hermite-Hadamard
type involving a three points divided difference (but not $f\left(\frac{x+y}{2}\right)$ 
$\frac{f(x)+f(y)}{2}$).

\begin{proposition}
\label{4}
Let $f:[x,y]\to\eR$ be a continuous and convex function and let $F$ be its antiderivative. Then $f$ and $F$ satisfy the following
inequality
\begin{equation}
\label{3}
f(x)\leq\frac{-3F(x)+4F\left(\frac{x+y}{2}\right)-F(y)}{y-x}.
\end{equation}
\end{proposition}
{\sc Proof.} For the sake of simplicity we shall assume that $x=0$ and $y=1.$
Then, according to Proposition \ref{l1}
$$\frac{-3F(x)+4F\left(\frac{x+y}{2}\right)-F(y)}{y-x}=\int_{0}^1fd\mu=\int_0^{\frac12}3f(t)dt-\int_{\frac12}^1f(t)dt.$$
Let
$$F_1(t):=\left\{\begin{array}{ll}
0&t<0\\
1&t\geq 0
\end{array}\right.$$
and let 
$$F_2(t):=\mu\bigl((-\infty,t]\bigr).$$
Then $(F_1,F_2)$ crosses only once (at $\frac13$) and we have 
$$\int_0^1tdF_1(t)=\int_0^1tdF_2(t)=0$$
 thus it suffices to use Ohlin lemma to get 
$$f(0)=\int_0^1fdF_1\leq\int_0^1fdF_2=-3F(0)+4F\left(\frac{1}{2}\right)-F(1).$$ 

As we can see, it is possible to prove some results for three point formulas but to get 
more interesting results we shall use longer expressions. For example, from numerical analysis, it is known that 
for $F'=f$ we have
\begin{equation}
\label{8}
f\left(\frac{x+y}{2}\right)\approx\frac{\frac13F(x)-\frac83F\left(\frac{3x+y}{4}\right)+\frac83F\left(\frac{x+3y}{4}
\right)-\frac13F(y)}{y-x}.
\end{equation}
It is natural to ask if the expression occurring at the right-hand side satisfies inequalities of the Hermite-Hadamard type. First we need some auxiliary results.

\begin{lemma} If any of the inequalities
\begin{equation}
\label{lhh}
f\left(\frac{x+y}{2}\right)\leq\frac{\sum_{i=1}^{4}a_iF(\alpha_i x+(1-\alpha_i)y)}{y-x}
\end{equation}
or 
\begin{equation}
\label{rhh}
\frac{\sum_{i=1}^{4}a_iF(\alpha_i x+(1-\alpha_i)y)}{y-x}\leq\frac{f(x)+f(y)}{2}
\end{equation}
is satisfied for all continuous and convex functions $f:[x,y]\to\eR$ (where $F'=f$) then 
\begin{equation}
\label{M1}
a_1(\alpha_2-\alpha_1)+(a_2+a_1)(\alpha_3-\alpha_2)+(a_3+a_2+a_1)(\alpha_4-\alpha_3)=1
\end{equation}
and 
\begin{equation}
\label{E1}
a_1(\alpha_2^2-\alpha_1^2)+(a_2+a_1)(\alpha_3^2-\alpha_2^2)+(a_3+a_2+a_1)(\alpha_4^2-\alpha_3^2)=1
\end{equation}
\end{lemma}
{\sc Proof.} Taking $x=0,y=1$ and, using Proposition \ref{l1}, we can see that 
$$\sum_{i=1}^{4}a_iF(1-\alpha_i)=\int_0^1fd\mu=-a_1\int_{1-\alpha_1}^{1-\alpha_2}f(x)dx+$$
$$-(a_1+a_2)\int_{1-\alpha_3}^{1-\alpha_2}f(x)dx-(a_1+a_2+a_3)\int_{1-\alpha_4}^{1-\alpha_3}f(x)dx.$$
Now we define $F_1,F_3$ and $F_4$ by the formulas \rf{F1},\rf{F3} and \rf{F4}, respectively. Then inequalities
\rf{lhh} and \rf{rhh} may be written in the form.
$$\int fdF_4\leq\int fdF_1$$
and
$$\int fdF_1\leq\int fdF_3.$$
This means that, if for example inequality \rf{lhh} is satisfied then we must have
$F_1(1)=F_4(1)=1$ which yields \rf{M1}. Further 
$$\int_{0}^1F_1(t)dt=\int_{0}^1F_4(t)dt=\frac12$$
which gives us \rf{E1}.

The above lemma gives only necessary conditions for inequalities \rf{lhh},\rf{rhh}. It is easy to see that the expression from \rf{8} satisfies \rf{M1} and \rf{E1}. To check if in concrete situations \rf{lhh} and 
\rf{rhh} are satisfied we have to use  Ohlin lemma or a more general Lemma \ref{mainl}.
To this end we formulate the following proposition.

\begin{proposition}
\label{ex1}
Let $\alpha_i, i=1,\dots,4$ satisfy $1=\alpha_{1}>\alpha_2>\alpha_3>\alpha_4=0$ let $a_i\in\eR$ be
such that $a_1+a_2+a_3+a_4=0$ and let equalities \rf{M1} and \rf{E1} be satisfied.
If $F_1$ is such that 
$$\frac{\sum_{i=1}^{4}a_iF(\alpha_{i}x+(1-\alpha_i)y)}{y-x}=\int_x^y fdF_1$$
and $F_2$ is the distribution function of a measure which is uniformly distributed in the interval $[x,y]$ then 
$(F_1,F_2)$ crosses exactly once.
\end{proposition}
{\sc Proof} From \rf{M1} we can see that $F_1(x)=F_2(x)=0$ and $F_1(y)=F_2(y)=1.$ Note that,
in view of Proposition \ref{l1} the graph of the restriction of $F_1$ to the interval $[x,y]$ consists of three segments. Therefore  $F_1$ and $F_2$ cannot have more than one crossing point. On the other hand if graphs $F_1$ and $F_2$ do not cross then 
$$\int_x^y tdF_1(t)\neq\int_x^y tdF_1(t)$$
i.e. \rf{E1} is not satisfied.

\begin{theorem}Let $\alpha_i, i=1,\dots,4$ satisfy $1=\alpha_{1}>\alpha_2>\alpha_3>\alpha_4=0,$ let $a_i\in\eR$ be
such that $a_1+a_2+a_3+a_4=0$ and let inequalities \rf{M1} and \rf{E1} be satisfied.
Let $F,f:[x,y]\to\eR$ be functions such that $f$ is continuous and convex and $F'=f.$ Then

(i) if $a_1>-1$ then
$$\frac{\sum_{i=1}^{4}a_iF(\alpha_ix+(1-\alpha_i)y)}{y-x}\leq\frac{1}{y-x}\int_x^yf(t)dt\leq\frac{f(x)+f(y)}{2}$$

(ii) if $a_1<-1$ then
$$f\left(\frac{x+y}{2}\right)\leq\frac{1}{y-x}\int_x^yf(t)dt\leq\frac{\sum_{i=1}^{4}a_iF(\alpha_ix+(1-\alpha_i)y)}{y-x}$$

(iii) if $a_1\in(-1,0]$ then
$$f\left(\frac{x+y}{2}\right)\leq\frac{\sum_{i=1}^{4}a_iF(\alpha_ix+(1-\alpha_i)y)}{y-x}\leq\frac{1}{y-x}\int_x^yf(t)dt$$

(iv) if $a_1<-1$ and $a_2+a_1\leq 0$ then
$$\frac{1}{y-x}\int_x^yf(t)dt\leq\frac{\sum_{i=1}^{4}a_iF(\alpha_ix+(1-\alpha_i)y)}{y-x}\leq\frac{f(x)+f(y)}{2}$$
\end{theorem}
{\sc Proof} We shall prove the first assertion. Other proofs are similar and will be omitted.  
It is easy to see that if inequalities which we consider are satisfied by every continuous and convex function
defined on the interval $[0,1]$ then they are true for every continuous and convex function on a given interval $[x,y].$ Therefore we assume that $x=0$ and $y=1.$ 
Let $F_1$ be such that \rf{Fint} is satisfied and let $F_2$ be the distribution function of a measure which is uniformly distributed in the interval $[0,1].$  From Proposition \ref{l1} and Remark 
\ref{R1} we can see that the graph of $F_1$ consists of three segments and, since $a_1>-1,$
the slope of the first segment is smaller than $1,$ i.e. $F_1$ lies below $F_2$ on some 
right-hand neighborhood of $x.$  In view of the Proposition \ref{ex1}, this means that the assumptions of  Ohlin lemma are satisfied and we get our result from this lemma. $\Box$

Now we shall present examples of inequalities which may be obtained from this theorem.

\begin{example}
\label{ex0}
Using (i), we can see that the inequality 
$${\frac13F(x)-\frac83F\left(\frac{3x+y}{4}\right)+\frac83F\left(\frac{x+3y}{4}
\right)-\frac13F(y)}\leq\frac{\int_x^yf(t)dt}{y-x}$$
is satisfied for every continuous and convex $f$ and its antiderivative $F.$
\end{example}

\begin{example}
\label{ex2}
Using (ii), we can see that the inequality 
$${-2F(x)+3F\left(\frac{2x+y}{3}\right)-3F\left(\frac{x+2y}{3}
\right)+2F(y)}\geq\frac{\int_x^yf(t)dt}{y-x}$$
is satisfied by every continuous and convex function $f$ and its antiderivative $F.$
\end{example}

\begin{example}Using (iii), we can see that the inequality 
$$\frac{\int_x^yf(t)dt}{y-x}\geq\frac{-\frac12F(x)-\frac32F\left(\frac{2x+y}{3}\right)+\frac32F\left(\frac{x+2y}{3}
\right)+\frac12F(y)}{y-x}\geq f\left(\frac{x+y}{2}\right)$$
is satisfied by every continuous and convex function $f$ and its antiderivative $F.$
\end{example}

\begin{example}Using (iv), we can see that the inequality 
$$\frac{\int_x^yf(t)dt}{y-x}\leq\frac{-\frac32F(x)+2F\left(\frac{3x+y}{4}\right)-2F\left(\frac{x+3y}{4}
\right)+\frac32F(y)}{y-x}\leq\frac{f(x)+f(y)}{2}$$
is satisfied by every continuous and convex function $f$ and its antiderivative $F.$
\end{example}

\begin{remark}
Inequalities occurring in Examples \ref{ex0} and \ref{ex2} may be written using function $F,$ exclusively. Thus we obtain, in fact, inequalities satisfied by $2-$convex function.
\end{remark}

In all cases considered in the above theorem we used only Ohlin lemma. Using Lemma \ref{mainl}, it is possible to 
obtain more subtle inequalities. However (for the sake of simplicity) in the next result we shall restrict our considerations to expressions of the simplified form. Note that 
the inequality between  $f\left(\frac{x+y}{2}\right)$ and 
expressions which we consider is a bit unexpected.

\begin{theorem}
\label{t2}
Let $\alpha\in\left(0,\frac12\right)$ let $a,b\in\eR$   
 and let inequalities \rf{M1} and \rf{E1} be satisfied.

(i) if $a>0$ then inequality 
$$f\left(\frac{x+y}{2}\right)\geq\frac{aF(x)+bF(\alpha x+(1-\alpha)y)-bF((1-\alpha) x+\alpha y)-aF(y)}{y-x}$$
is satisfied by every continuous and convex $f$ and its antiderivative $F$ if and only if
\begin{equation}
\label{A1A2}
(1-\alpha)^2\frac{ab}{a+b}>\frac12-(1-\alpha)\frac{b}{a+b},
\end{equation}

(ii) if $a<-1$ and $a_1+a_2>0$ then inequality 
\begin{equation}
\label{>}
\frac{aF(x)+bF(\alpha x+(1-\alpha)y)-bF((1-\alpha) x+\alpha y)-aF(y)}{y-x}\leq \frac{f(x)+f(y)}{2}
\end{equation}
is satisfied by every continuous and convex $f$ and its antiderivative $F$ if and only if
$$-\frac{1}{4a}>\left(-a(1-\alpha)-\frac12\right)\left(\frac12+\frac{1}{2a}\right).$$
\end{theorem}
{\sc Proof} We shall prove the assertion (i). The proof of (ii) is similar and will be omitted. 
Similarly as before we assume without loss of generality that $x=0,y=1$ and let $F_1$ be such that 
$$aF(0)+bF(1-\alpha)-bF(\alpha)+aF(1)=\int_0^1fdF_1$$
further let $F_4$ be given by \rf{F4}. Then it is easy to see that $(F_1,F_4)$ crosses three times:
at $\frac{(1-\alpha)b}{a+b},$ at $\frac12$ and at $\frac{a+\alpha b}{a+b},$

We are going to use Lemma \ref{maini}. Since from \rf{E1} we know that
$$A_0+A_1+A_2+A_3=0$$ 
it suffices to check that 
$A_0\geq A_1$ if and only if inequality \rf{A1A2} is satisfied.
Since $F_4(x)=0,$ for  $x\in\left(0,\frac12\right),$ we get 
$$A_0=-\int_{0}^{\frac{(1-\alpha)b}{a+b}}F_1(t)dt$$
and 
$$A_1=\int_{\frac{(1-\alpha)b}{a+b}}^\frac12F_1(t)dt$$
which yields our assertion.

\begin{example}Neither inequality
\begin{equation}
\label{*}
f\left(\frac{x+y}{2}\right)\leq\frac{\frac13F(x)-\frac83F\left(\frac{3x+y}{4}\right)+\frac83F\left(\frac{x+3y}{4}
\right)-\frac13F(y)}{y-x}
\end{equation}
 nor 
 \begin{equation}
 \label{**}
f\left(\frac{x+y}{2}\right)\geq\frac{\frac13F(x)-\frac83F\left(\frac{3x+y}{4}\right)+\frac83F\left(\frac{x+3y}{4}
\right)-\frac13F(y)}{y-x}
\end{equation}
is satisfied for all continuous and convex $f:[x,y]\to\eR.$ Indeed, if $F_1$ is such that  
$$\int_x^yf(t)dF_1(t)=\frac{\frac13F(x)-\frac83F\left(\frac{3x+y}{4}\right)+\frac83F\left(\frac{x+3y}{4}
\right)-\frac13F(y)}{y-x}$$
then 
$$\int_x^{\frac{3x+y}{4}}F_1(t)dt<\int_x^{\frac{3x+y}{4}}F_4(t)dt$$
thus inequality \rf{*} cannot be satisfied. On the other hand, the coefficients and nodes 
of the expression considered do not satisfy \rf{A1A2}. Therefore \rf{**} is also not satisfied 
for all continuous and convex $f:[x,y]\to\eR.$
\end{example}

\begin{example}Using Theorem \ref{t2} assertion (i), we can see that the inequality 
$$\frac{2F(x)-3F\left(\frac{3x+y}{4}\right)+3F\left(\frac{x+3y}{4}
\right)-2F(y)}{y-x}\leq f\left(\frac{x+y}{2}\right)$$
is satisfied for every continuous and convex $f$ and its antiderivative $F.$
\end{example}

\begin{example}Using Theorem \ref{2} assertion (ii), we can see that the inequality 
$$\frac{-2F(x)+3F\left(\frac{2x+y}{3}\right)-3F\left(\frac{x+2y}{3}
\right)+2F(y)}{y-x}\leq \frac{f(x)+f(y)}{2}$$
is satisfied for every continuous and convex $f$ and its antiderivative $F.$
\end{example}

\begin{remark} As it is known from the paper \cite{BesPal}, if a continuous function
satisfies inequalities  
of the type which we have considered then such function must be convex.

Therefore inequalities obtained in this paper characterize convex functions 
(in the class of continuous functions).
\end{remark}

%\begin{acknowledgements}
%\end{acknowledgements}



\end{document}